\documentclass[12pt,twoside]{article}

\usepackage{amsmath,amssymb,amsthm,amscd,graphicx,color,accents}

\usepackage[colorlinks=true,linkcolor=blue,citecolor=blue,urlcolor=blue]{hyperref}
\usepackage[round]{natbib}

\numberwithin{equation}{section}

\theoremstyle{plain}
\newtheorem{theorem}{Theorem}[section]
\newtheorem{proposition}[theorem]{Proposition}

\newtheorem{corollary}[theorem]{Corollary}
\newtheorem{remark}[theorem]{Remark}

\newtheorem{example}[theorem]{Example}

\newcommand\wtilde{\widetilde}
\newcommand{\Pro}{\noindent\textit{Proof.}\ \ }

\def\e{{\rm{e}}}

\def\bB{{\boldsymbol{B}}}
\def\bS{{\boldsymbol{S}}}

\def\bLambda{{\boldsymbol{\Lambda}}}
\def\bOmega{{\boldsymbol{\Omega}}}
\def\dd{\,{\rm{d}}}

\def\E{\mathbb{E}}

\def\N{\mathbb{N}}
\def\P{\mathbf{P}}
\def\R{\mathbb{R}}

\begin{document}

\title{\Large\textbf{Loading Monotonicity of Weighted Premiums, and Total Positivity Properties of Weight Functions}}

\author{{Donald Richards}{\hskip1pt}\thanks{
Department of Statistics, Pennsylvania State University, University Park, PA 16802, U.S.A.  E-mail address: \href{mailto:richards@stat.psu.edu}{richards@stat.psu.edu}.
\endgraf
\ $^\dag$Laboratory for Information and Decision Systems and Institute for Data, Systems and Society, Massachusetts Institute of Technology, Cambridge, MA 02139, U.S.A. \ E-mail address: \href{mailto:cuhler@mit.edu}{cuhler@mit.edu}.
\endgraf
\endgraf
\endgraf
\endgraf
}\ \ and {Caroline Uhler}{\hskip1pt}$^\dag$
\endgraf
}

\date{\today}

\maketitle

\normalsize

\begin{abstract}
We consider the construction of insurance premiums that are monotonically increasing with respect to a loading parameter.  By introducing weight functions that are totally positive of higher order, we derive higher monotonicity properties of generalized weighted premiums; in particular, we deduce for weight functions that are totally positive of order three a monotonicity property of the variance-to-mean ratio, or index of dispersion, of the loss variable.  We derive the higher order total positivity properties of some ratios that arise in actuarial and insurance analysis of combined risks.  Further, we examine seven classes of weight functions that have appeared in the literature and we ascertain the higher order total positivity properties of those functions.
\end{abstract}

\section{Introduction}
\label{introduction}

Consider the problem of estimating the premiums that an insurance operation is to charge its clients in order to underwrite their risks.  On the one hand, the insurer is limited by competition as to how much it may charge to underwrite a given risk.  On the other hand the insurer, so as to remain solvent, necessarily must charge premiums that are suitably large in order to cover its insured risks and its operating expenses.  

To formulate this problem probabilistically, suppose that we have a probability triplet $(\bOmega,\mathcal{A},\P)$, consisting of a sample space $\bOmega$, a sigma-algebra $\mathcal{A}$ of subsets of $\bOmega$, and a probability measure $\P(\cdot)$ on $\mathcal{A}$.  We suppose that there corresponds to $\P(\cdot)$ a random variable $X : \bOmega \to \R_+$, called the {\it loss variable}, that arises when the insurer underwrites a randomly chosen risk.  We assume that $X$ is nonnegative and refer to its mean, $\E[X] = \int_\bOmega X(\omega) \P(\dd\omega)$, as the {\it net premium}.  Noting that the net premium $\E[X]$ will cover only the average insured risk, the insurer, in order to remain profitable, necessarily must charge an amount $H[X]$ that is {\it loaded}, meaning that $H[X] \ge \E[X]$.  

We refer to \cite{furman_zitikis08a} for an extensive review of the construction of loaded premiums, and to \cite{sendov_etal} for further results on the same topic.  We also cite \cite{jones_zitikis} and \cite{furman_zitikis08b, furman_zitikis08c, furman_zitikis09} for motivating accounts and related analyses of actuarial and insurance problems that involve loading parameters.  

A method for constructing weighted premiums begins with the insurer choosing a nonnegative {\it weight function} $w(\lambda,x)$ that depends on a {\it loading parameter} $\lambda > 0$.  We are especially interested in the weighted premium, 
$$%\begin{equation}
%\label{Hfunction}
H[\lambda,X] = \frac{\E[X \, w(\lambda,X)]}{\E[w(\lambda,X)]},
$$%\end{equation}
where it is assumed that, for each $\lambda > 0$, the function $x \mapsto w(\lambda,x)$ is Borel-measurable.  Suppose that the weight function $(\lambda,x) \mapsto w(\lambda,x)$ is {\it totally positive of order} $2$, i.e., 
$$%\begin{equation}
%\label{tp2}
w(\lambda_1,x_1) \, w(\lambda_2,x_2) \ge w(\lambda_1,x_2) \, w(\lambda_2,x_1),
$$%\end{equation}
whenever $\lambda_1 > \lambda_2$ and $x_1 > x_2$; then \citet[Theorem 2.1]{sendov_etal} proved that $H[\lambda,X]$ is non-decreasing in $\lambda$, thereby relating the study of weighted premiums with the theory of total positivity.  The implication for insurance pricing is that if $w(\lambda,x)$ is totally positive of order two then a riskier venture, with risk represented by the parameter $\lambda$, will not be assigned a lower weighted net premium.  We refer to \cite{furman_zitikis08a} who introduced the concept of a {\it weighted premium} in research on the construction of insurance premiums.  

Noting the general theory of total positivity \citep{karlin}, we wish to determine the behavior of the weighted premium $H[\lambda,X]$ for weights $w(\lambda,x)$ that are totally positive of order higher than two.  In this paper, we study $H[\lambda,X]$, and some of its generalizations, when the weight function $w(\lambda,x)$ is totally positive of any given order.  

Our results may be described as follows.  In Section \ref{totalpositivity}, we introduce the theory of total positivity, providing a self-contained introduction to results needed in the sequel.  

We consider in Section \ref{generalizedmonotonicity} classes of {\it generalized weighted premiums}, as defined by \cite{furman_zitikis09}, extending $H[\lambda,X]$.  We establish monotonicity properties of the generalized weighted premiums, recovering as a special case the previously cited result of \citet[Theorem 2.1]{sendov_etal}, and we deduce for weight functions that are totally positive of order three a monotonicity property of the variance-to-mean ratio (or index of dispersion) of the loss variable $X$.  In Section \ref{actuarialratios}, we derive some total positivity properties of $\mathcal{R}_c$ and $\mathcal{C}_c$, two actuarial ratios that were defined and studied by \cite{furman_zitikis08b,furman_zitikis08c} in the analysis of combined risks.  

In Section \ref{sevenweightfunctions}, we consider seven classes of weight functions treated previously by \cite{sendov_etal}.  We ascertain the higher order total positivity properties of these weight functions, proving that five of them are strictly totally positive of order infinity, one is totally positive of order infinity, and one is not totally positive of order three.  

Finally, in Section \ref{sec:conclusions}, we summarize with concluding remarks on the implications of working with weighted premiums that are totally positive of higher order.

\section{Total positivity}
\label{totalpositivity}

We begin by recalling from \cite{karlin} the concepts of total positivity, strict total positivity, and sign regularity.  

For $k \in \N$, a weight function $w:\R^2 \to \R$ is {\it totally positive of order} $k$, denoted TP$_k$, if for all $\lambda_1 > \cdots > \lambda_k$, $x_1 > \cdots > x_k$, and for all $r=1,\ldots,k$, the $r \times r$ determinant, 
$$
\det\big(w(\lambda_i,x_j)\big) := 
\left|
\begin{matrix}
w(\lambda_1,x_1) & \cdots & w(\lambda_1,x_r) \\
\vdots & \cdots & \vdots \\
w(\lambda_r,x_1) & \cdots & w(\lambda_r,x_r)
\end{matrix}
\right| \ge 0.
$$
The function $w(\lambda,x)$ is {\it totally positive of order infinity}, denoted TP$_\infty$, if $w(x,\lambda)$ is TP$_k$ for all $k \ge 1$.  Similarly, $w(\lambda,x)$ is {\it strictly totally positive of order} $k$, denoted STP$_k$ if the $r \times r$ determinant $\det\big(w(\lambda_i,x_j)\big)$ is strictly positive for all $\lambda_1 > \cdots > \lambda_k$, $x_1 > \cdots > x_k$, and all $r=1,\ldots,k$.  Further, $w(\lambda,x)$ is {\it strictly totally positive of order infinity}, denoted STP$_\infty$, if $w(\lambda,x)$ is STP$_k$ for all $k \ge 1$.  

The function $w(\lambda,x)$ is said to be {\it reverse-rule of order} $k$, denoted RR$_k$, if for all $\lambda_1 > \cdots > \lambda_k$ and $x_1 > \cdots > x_k$, $(-1)^{r(r-1)/2} \det\big(w(\lambda_i,x_j)\big)$ is nonnegative for all $r=1,\ldots,k$; if this holds for all $k \ge 1$ then $w(\lambda,x)$ is called {\it reverse-rule of order infinity}, denoted RR$_\infty$.  If $(-1)^{r(r-1)/2} \det\big(w(\lambda_i,x_j)\big)$ is strictly positive for all $\lambda_1 > \cdots > \lambda_k$, $x_1 > \cdots > x_k$, $r=1,\ldots,k$ then $w(\lambda,x)$ is said to be {\it strictly reverse-rule of order} $k$ (SRR$_k$); and $w(\lambda,x)$ is called {\it strictly reverse-rule of order infinity} (SRR$_\infty$) if it is SRR$_k$ for all $k \ge 1$.

%Examples of TP$_\infty$ or STP$_\infty$ functions: $\e^{\lambda x}$, $x^\lambda$, $\e^{-x/\lambda}$

Throughout the remainder of the paper, we will assume that all integrals or sums converge absolutely.  Whenever it is necessary to provide explicit conditions under which such convergence holds then we will provide the details.  

The {\it Binet-Cauchy formula} often is stated in terms of calculating the minors of a matrix product, $AB$, from the minors of $A$ and $B$ \cite[p.~1]{karlin}.  We will need a continuous and a discrete generalization of this formula:  Let $\nu$ be a Borel-finite measure on a totally ordered measure space $\mathfrak{X}$.  Also, for $r \in \N$, let $\phi_1,\ldots,\phi_r$ and $\psi_1,\ldots,\psi_r$ be complex-valued functions on $\mathfrak{X}$.  The Binet-Cauchy formula is that the $r \times r$ determinant with $(i,j)$th entry $\int_{\mathfrak{X}} \phi_i(x) \psi_j(x) \dd\nu(x)$ satisfies the identity 
\begin{equation}
\label{binetcauchy_continuous}
\det\bigg(\int_{\mathfrak{X}} \phi_i(x) \psi_j(x) \dd\nu(x)\bigg) = \mathop{\idotsint}_{x_1 > \cdots > x_r} \det\big(\phi_i(x_j)\big) \det\big(\psi_i(x_j)\big) \prod_{j=1}^r \dd\nu(x_j).
\end{equation}
%whenever all integrals appearing in (\ref{binetcauchy_continuous}) converge absolutely.  

For the case in which $\mathfrak{X} = \N_0$, the set of nonnegative integers, and $\nu$ is a discrete measure on $\N_0$ with weights $\nu(m)$, $m=0,1,2,\ldots$, the Binet-Cauchy formula is the statement that 
\begin{equation}
\label{binetcauchy_discrete}
\det\bigg(\sum_{m=0}^\infty \phi_i(m) \psi_j(m) \nu(m)\bigg) = \mathop{\sum}_{m_1 > \cdots > m_r \ge 0} \det\big(\phi_i(m_j)\big) \det\big(\psi_i(m_j)\big) \prod_{j=1}^r \nu(m_j).
\end{equation}
%whenever all sums in (\ref{binetcauchy_discrete}) converge absolutely.  

The continuous version of the {\em Basic Composition Formula} is that if the weight functions $w_1(\lambda,x)$ and $w_2(\lambda,x)$ are TP$_k$ on $\R^2$, and if $\nu$ is a sigma-finite measure on $\R$, then the weight function 
\begin{equation}
\label{bcf_continuous}
w(\lambda,x) = \int_\R w_1(\lambda,t) w_2(t,x) \dd\nu(t)
\end{equation}
also is TP$_k$ on $\R^2$.  

The discrete version of the Basic Composition Formula, analogous to (\ref{binetcauchy_discrete}), is that if $w_1(\lambda,x)$ and $w_2(\lambda,x)$ are TP$_k$ on $\N_0 \times \N_0$, and $\nu$ is a discrete measure on $\N_0$ with nonnegative weights $\nu(m)$, $m=0,1,2,\ldots$, then the function 
\begin{equation}
\label{bcf_discrete}
w(\lambda,x) = \sum_{m=0}^\infty w_1(\lambda,m) w_2(m,x) \nu(m)
\end{equation}
also is TP$_k$ on $\N_0 \times \N_0$.  

We remark that a crucial difference between total positivity of order two and total positivity of higher orders is that if a positive function $w(\lambda,x)$ is TP$_2$ then the function $1/w(\lambda,x)$ is RR$_2$. However this result does not generally extend to TP$_k$ functions for $k > 2$. This explains why some of the weight functions considered in Section \ref{sevenweightfunctions} have relatively straightforward TP$_2$ or RR$_2$ properties, while their higher-order total positivity properties are more difficult to establish.  We refer to \citet[Eq.~(1.5)]{carlson_gustafson} for further remarks on this point.

%Alternative terminology: total positivity vs. log-supermodularity. 

\section{Monotonicity properties of generalized weighted premiums}
\label{generalizedmonotonicity}

Let $X$ be a nonnegative random variable with probability density function $g$.  \citet{sendov_etal} derived a monotonicity property of the weighted premium function, 
$$
H[\lambda,X] = \frac{\E[w(\lambda,X) \, X]}{\E[w(\lambda,X)]},
$$
where the expectations are taken with respect to the distribution of $X$.  \citet[Theorem 2.1]{sendov_etal} proved that if $w(\lambda,x)$ is TP$_2$ then the function $\lambda \mapsto H[\lambda,X]$ is non-decreasing.  We shall generalize this property in two ways.  Following \cite{furman_zitikis09}, we consider for a {\it utility function} $f: \R_+ \to \R_+$, the {\it generalized weighted premium},  
$$
H[\lambda,f(X)] = \frac{\E [w(\lambda,X) f(X)]}{\E[w(\lambda,X)]},
$$
whenever these expectations exist.  Let $Y$ be the random variable that has the weighted probability density function, 
\begin{equation}
\label{weighteddistn}
\frac{w(\lambda,y)}{\E[w(\lambda,X)]} \, g(y)
\end{equation}
$y \ge 0$, where, as defined earlier, $g$ is the density function of $X$.  Then $H[\lambda,f(X)]$ can also be viewed as the expectation $\E_Y f(Y)$, where the expectation is with respect to the distribution of $Y$.  We will establish the monotonicity of $H[\lambda,f(X)]$ for the case in which $f$ is monotonically increasing.  

Second, for the case in which the weight function $w(\lambda,x)$ is TP$_k$ or STP$_k$, we obtain generalizations of the monotonicity property arising from the case in which $k=2$.

\begin{theorem}
\label{generalizedloadingmonotonicity}
Suppose that the weight function $w(\lambda,x)$  is TP$_k$, $f: \R_+ \to \R_+$ is a non-decreasing function, and $\lambda_1 > \cdots > \lambda_k$.  Then, all minors of the $k \times k$ determinant 
\begin{equation}
\label{higherHfunctions}
\det\Big(H\big[\lambda_i,\big(f(X)\big)^{k-j}\,\big]\Big)
\end{equation}
are nonnegative.  Further, if $w(\lambda,x)$ is STP$_k$ and the set of points of increase of $f$ contains an open set then all minors of the matrix (\ref{higherHfunctions}) are positive.  
%In particular, $\lambda \mapsto H[\lambda,f(X)]$ is strictly increasing.  
\end{theorem}

\Pro
Consider the $m \times m$ minor of (\ref{higherHfunctions}) corresponding to rows $r_1,\ldots,r_m$ and columns $c_1,\ldots,c_m$, where $r_1 < \cdots < r_m$ and $c_1 < \cdots < c_m$.  Applying the Binet-Cauchy formula (\ref{binetcauchy_continuous}) with $\mathfrak{X} = \R$, $\phi_i(x) = w(\lambda_{c_i},x)$ and $\psi_i(x) = \big(f(x)\big)^{k-r_i}$, $i=1,\ldots,m$, and $\dd \nu(x) = g(x) \dd x$, we obtain 
\begin{multline}
\label{binetcauchyhigherH}
\det\big(\E[w(\lambda_{c_i},X) \big(f(X)\big)^{k-r_j}]\big) \\
= \idotsint\limits_{x_1 > \cdots > x_m} \det\big(w(\lambda_{c_i},x_j)\big) \det\Big(\big(f(x_j)\big)^{k-r_i}\Big) \prod_{j=1}^m g(x_j) \dd x_j.
\end{multline}
Since $w(\lambda,x)$ is TP$_r$ and $\lambda_{c_1} > \cdots > \lambda_{c_m}$ then $\det\big(w(\lambda_{c_i},x_j)\big)$ is nonnegative on the orthant $\{(x_1,\ldots,x_m): x_1 > \cdots > x_m\}$.  

As for the second determinant in the integrand in (\ref{binetcauchyhigherH}), let $\theta_i = k-r_i-m+i$, $i=1,\ldots,m$, and set $\theta = (\theta_1,\ldots,\theta_m)$.  Since $1 \le r_1 < \cdots < r_m \le k$ then $k - m \ge \theta_1 \ge \cdots \ge \theta_m \ge 0$.  The determinant 
$$
\det\big(t_j^{k-r_i}\big) \equiv \det\big(t_j^{\theta_i + m - i}\big)
$$
is well-known; see \citet[p.~40]{macdonald}.  In particular, this determinant is divisible by the product $\prod_{1 \le i < j \le m} (t_i - t_j)$, and the ratio of these two polynomials defines the {\it Schur function}, 
\begin{equation}
\label{schurdef}
\chi_\theta(t_1,\ldots,t_m) = \frac{\det\big(t_j^{\theta_i + m - i}\big)}{\prod_{1 \le i < j \le m} (t_i - t_j)}.
\end{equation}

It is straightforward to verify that $\chi_\theta(t_1,\ldots,t_m)$ is a homogeneous polynomial of degree $\theta_1+\cdots+\theta_m$.  It is also well-known that the coefficients appearing in the monomial expansion of $\chi_\theta(t_1,\ldots,t_m)$ are nonnegative integers \citep[p.~75]{macdonald}.  Therefore, $\chi_\theta(t_1,\ldots,t_m) > 0$ for $t_1,\ldots,t_m > 0$.  Writing (\ref{schurdef}) in the form 
$$
\det\big(t_j^{k-r_i}\big) = \prod_{1 \le i < j \le m} (t_i - t_j) \cdot \chi_\theta(t_1,\ldots,t_m),
$$
it follows that $\det\big(t_j^{k-r_i}\big) > 0$ for all $t_1,\ldots,t_m > 0$.  Consequently, by substituting $t_i = f(x_i)$, we obtain 
\begin{equation}
\label{schur}
\det\Big(\big(f(x_j)\big)^{k-r_i}\Big) = \prod_{1 \le i < j \le m} \big(f(x_i) - f(x_j)\big) \cdot \chi_\theta\big(f(x_1),\ldots,f(x_m)\big),
\end{equation}
and since $f$ is increasing then it follows that the determinant in (\ref{schur}) is nonnegative for $x_1 > \cdots > x_m$.  

Therefore, the integrand in (\ref{binetcauchyhigherH}) is nonnegative for $\lambda_{c_1} > \cdots > \lambda_{c_m}$ and $x_1 > \cdots > x_m$, so it follows that $\det\big(\E\big[w(\lambda_{c_i},X) \big(f(X)\big)^{k-r_j}\big]\big) \ge 0$.  Since $m$, $r_1,\ldots,r_m$ and $c_1,\ldots,c_m$ were chosen arbitrarily then we deduce that all minors of the $k \times k$ determinant $\det\big(\E\big[w(\lambda_i,X) \big(f(X)\big)^{k-j}\big]\big)$ are nonnegative.  

If $w(\lambda,x)$ is STP$_k$ then $\det\big(w(\lambda_{c_i},x_j)\big) > 0$ for all $\lambda_{c_1} > \cdots > \lambda_{c_m}$ and $x_1 > \cdots > x_m$.  If also the set of points of increase of $f$ contains an open set then the determinant (\ref{schur}) is positive on an open set in the orthant $\{(x_1,\ldots,x_m): x_1 > \cdots > x_m\}$.  Then, the integrand in (\ref{binetcauchyhigherH}) is positive on an open set, so it follows that $\det\big(\E\big[w(\lambda_{c_i},X) \big(f(X)\big)^{k-r_j}\big]\big) > 0$.  

For $j=1,\ldots,r$, we divide by $\E\big[w(\lambda_{c_i},X)\big]$ the $j$th column of the determinant $\det\big(\E\big[w(\lambda_{c_i},X) \big(f(X)\big)^{k-r_j}\big]\big)$.  Since 
$$
\frac{\E\big[w(\lambda_{c_i},X) \big(f(X)\big)^{k-r_j}\big]}{\E[w(\lambda_{c_i},X)]} = H[\lambda_{c_i},\big(f(X)\big)^{k-r_j}]
$$
then we find that $\det\big(H[\lambda_{c_i},\big(f(X)\big)^{k-r_j}]\big) \ge 0$ for all $\lambda_{c_1} > \cdots > \lambda_{c_m}$.  As before, it follows that all minors of $\det\big(H[\lambda_i,\big(f(X)\big)^{k-j}]\big)$ are nonnegative.  

Finally, for the case in which $w(\lambda,x)$ is STP$_k$ and $f$ is strictly increasing on an open set, we deduce analogously that $\det\big(H[\lambda_{c_i},\big(f(X)\big)^{k-r_j}]\big) > 0$ for all $\lambda_{c_1} > \cdots > \lambda_{c_m}$.  Therefore, all minors of $\det\big(H[\lambda_i,\big(f(X)\big)^{k-j}]\big)$ are positive.  
$\qed$

\smallskip

\begin{remark} 
\label{remarkhighermoments}
{\rm 
(1) Consider the case in which $k=2$.  As $H[\lambda,1] \equiv 1$, Theorem \ref{generalizedloadingmonotonicity} provides that if $f$ is increasing, $w(\lambda,x)$ is TP$_2$, and if $\lambda_1 > \lambda_2$ then
$$
\left|
\begin{matrix}
H[\lambda_1,f(X)] & H[\lambda_1,1] \\
H[\lambda_2,f(X)] & H[\lambda_2,1]
\end{matrix}
\right| = 
H[\lambda_1,f(X)] - H[\lambda_2,f(X)] \ge 0;
$$
that is, the function $\lambda \mapsto H[\lambda,f(X)]$ is non-decreasing.  For the case in which $f(x) = x$, we recover the result of \citet[Theorem 2.1]{sendov_etal}.  

Let $\mu_\lambda$ denote $H[\lambda,f(X)]$; equivalently, $\mu_\lambda$ is the mean of $f(Y)$ with respect to the weighted distribution (\ref{weighteddistn}).  Then the hypothesis that $w(\lambda,x)$ is TP$_2$ leads to the conclusion that $\mu_\lambda$ is increasing in $\lambda$.  

(2) Suppose that $k = 3$; then Theorem \ref{generalizedloadingmonotonicity} provides that if $f$ is increasing, $w(\lambda,x)$ is TP$_3$, and if $\lambda_1 > \lambda_2 > \lambda_3$ then all minors of the determinant 
$$
\left|
\begin{matrix}
H[\lambda_1,(f(X))^2] & H[\lambda_1,f(X)] & 1 \\
H[\lambda_2,(f(X))^2] & H[\lambda_2,f(X)] & 1 \\
H[\lambda_3,(f(X))^2] & H[\lambda_3,f(X)] & 1
\end{matrix}
\right|
$$
are nonnegative.  In particular, the $2 \times 2$ minor, 
\begin{multline}
\label{minortwobytwo}
\left|
\begin{matrix}
H[\lambda_1,(f(X))^2] & H[\lambda_1,f(X)] \\
H[\lambda_2,(f(X))^2] & H[\lambda_2,f(X)]
\end{matrix}
\right| \\
= H[\lambda_1,(f(X))^2] H[\lambda_2,f(X)] - H[\lambda_1,f(X)] H[\lambda_2,(f(X))^2]
\end{multline}
is nonnegative.  

Suppose that the weight function $w(\lambda,x)$ is differentiable in $\lambda$.  Also, suppose that its partial derivative, $\partial w(\lambda,x)/\partial\lambda$, is integrable and that 
$$
H_1[\lambda,f(X)]: = \frac{\partial}{\partial\lambda} H[\lambda,f(X)]
$$
exists.  Dividing (\ref{minortwobytwo}) by $\lambda_1 - \lambda_2$ and then letting $\lambda_1,\lambda_2 \to \lambda$, we obtain 
\begin{align*}
0 & \le \lim_{\lambda_1,\lambda_2 \to \lambda} \frac{H[\lambda_1,(f(X))^2] H[\lambda,f(X)] - H[\lambda_1,f(X)] H[\lambda,(f(X))^2]}{\lambda_1-\lambda_2} \\
& = H_1[\lambda,(f(X))^2] H[\lambda,f(X)] - H_1[\lambda,f(X)] H[\lambda,(f(X))^2] \\
& = H[\lambda,f(X)] H[\lambda,(f(X))^2] \, \frac{\partial}{\partial\lambda} \log \frac{H[\lambda,(f(X))^2]}{H[\lambda,f(X)]};
\end{align*}
equivalently, 
$$
\frac{\partial}{\partial\lambda} \log \frac{H[\lambda,(f(X))^2]}{H[\lambda,f(X)]} \ge 0.
$$
Hence, if $f:\R_+ \to \R_+$ is increasing and $w(\lambda,x)$ is TP$_3$ then the function $\lambda \mapsto H[\lambda,(f(X))^2]/H[\lambda,f(X)]$ is increasing.  

To interpret this result within a statistical context, define 
$$
\sigma^2_\lambda := H[\lambda,(f(X))^2] - (H[\lambda,f(X)])^2,
$$
representing the variance of $f(Y)$ with respect to the weighted distribution (\ref{weighteddistn}).  Then, 
$$
\frac{H[\lambda,(f(X))^2]}{H[\lambda,f(X)]} = \frac{\sigma^2_\lambda + \mu_\lambda^2}{\mu_\lambda} = \frac{\sigma^2_\lambda}{\mu_\lambda} +  \mu_\lambda.
$$
The ratio $\sigma^2_\lambda/\mu_\lambda$ is known classically as the {\it variance-to-mean ratio} or {\it index of dispersion} \cite[p.~72]{cox_lewis}, and we denote it by VMR{\hskip1pt}$_\lambda$.  The variance-to-mean ratio is a normalized measure of the extent to which the possible values of $X$ are dispersed, so that smaller values of VMR{\hskip1pt}$_\lambda$ correspond to more concentrated clustering of the possible values of $X$.  Thus, the assumption that $w(\lambda,x)$ is TP$_3$ leads to the conclusions that $\mu_\lambda$ and VMR{\hskip1pt}$_\lambda + \mu_\lambda$ are increasing functions of $\lambda$.  In the context of premium calculations, the variance-to-mean ratio was studied earlier by \cite{heilmann}.  

We now see that as $k$, the order of total positivity of the weight function $w(\lambda,x)$, increases, we are able to deduce correspondingly more intricate aspects of the monotonicity properties of $H[\lambda,f(X)]$ as a function of $\lambda$.  An implication of the above remark is that if an insurer expects greater variance-to-mean ratios for increasing values of the loading parameter $\lambda$ then it would be advisable to calculate premiums using weight functions that are STP$_k$ with $k \ge 3$.  
}\end{remark}

For $k=2$, a consequence of the proof of Theorem \ref{generalizedloadingmonotonicity} is that it provides in (\ref{binetcauchyhigherH}) an explicit representation for the difference $H[\lambda_1,f(X)] - H[\lambda_2,f(X)]$ as the integral of a nonnegative function, viz., 
\begin{multline}
\label{prelipschitz}
\E[w(\lambda_1,X)] \, \E[w(\lambda_2,X)] \, \big(H[\lambda_1,f(X)] - H[\lambda_2,f(X)]\big) \\
= \ \ \mathop{\int\int}\limits_{x_1 > x_2} \big(f(x_1) - f(x_2)\big) \, 
\left|\begin{matrix}
w(\lambda_1,x_1) & w(\lambda_1,x_2) \\
w(\lambda_2,x_1) & w(\lambda_2,x_2)
\end{matrix}\right|
g(x_1) g(x_2) \dd x_1 \dd x_2.
\end{multline}
Hence, the nonnegativity of the difference, $H[\lambda_1,f(X)] - H[\lambda_2,f(X)]$, is obtained immediately.  Further, as the following result shows, the integral representation (\ref{prelipschitz}) combined with an estimate on the variation of $f(x)$ leads to an upper bound on $H[\lambda_1,f(X)] - H[\lambda_2,f(X)]$.  

\smallskip

\begin{corollary}
\label{lipschitz}
Suppose that $(\lambda,x) \mapsto w(\lambda,x)$ is TP$_2$ and $f:\R_+ \to \R_+$ is a non-decreasing function that satisfies a uniform Lipschitz condition of order $1$, viz., 
$$
|f(x_1) - f(x_2)| \le |x_1 - x_2|
$$ 
for all $x_1$ and $x_2$.  Then, for all $\lambda_1 > \lambda_2$, 
\begin{equation}
\label{lipschitz1}
H[\lambda_1,f(X)] - H[\lambda_2,f(X)] \le H[\lambda_1,X] - H[\lambda_2,X].
\end{equation}
\end{corollary}

\Pro
Applying to (\ref{prelipschitz}) the Lipschitz condition on $f$, we obtain 
\begin{align}
\E[w(\lambda_1&,X)] \, \E[w(\lambda_2,X)] \, \big(H[\lambda_1,f(X)] - H[\lambda_2,f(X)]\big) 
\label{lipschitz2} \\
& \le \ \ \mathop{\int\int}\limits_{x_1 > x_2} (x_1 - x_2) 
\left|\begin{matrix}
w(\lambda_1,x_1) & w(\lambda_1,x_2) \\
w(\lambda_2,x_1) & w(\lambda_2,x_2)
\end{matrix}\right|
g(x_1) g(x_2) \dd x_1 \dd x_2 \nonumber \\
& \equiv \ \ \mathop{\int\int}\limits_{x_1 > x_2} 
\det\big(x_i^{2-j}\big) \cdot
\det\big(w(\lambda_i,x_j)\big) 
g(x_1) g(x_2) \dd x_1 \dd x_2.
\label{lipschitz3}
\end{align}
Applying the Binet-Cauchy formula (\ref{binetcauchy_continuous}), we deduce that (\ref{lipschitz3}) equals 
\begin{align}
\label{lipschitz4}
\det\big(\E[w(\lambda_i,X) X^{2-j}]\big) &= \det\big(\E[w(\lambda_i,X)] H[\lambda_i,X^{2-j}]\big) \nonumber \\
&= \E[w(\lambda_1,X)] \, \E[w(\lambda_2,X)] \, \big(H[\lambda_1,X] - H[\lambda_2,X]\big).
\end{align}
On comparing (\ref{lipschitz2}) and (\ref{lipschitz4}), and clearing the common terms on each side of that inequality, we obtain (\ref{lipschitz1}).
$\qed$

\smallskip

\begin{remark}
{\rm \citet[Section 4]{furman_zitikis08a} provide examples of utility functions that are of the form $f(x) = \int_0^x h(x) \dd x$, where $h(x) \ge 0$ for all $x$.  Suppose that $h$ is uniformly bounded with $h(x) \le 1$ for all $x$; examples of such $h$ are the cumulative distribution functions of nonnegative random variables.  Then for $x_1, x_2 \in \R$, it follows from the triangle inequality that 
$$
|f(x_1) - f(x_2)| 
%= \left|\int_0^{x_1} h(x) \dd x - \int_0^{x_2} h(x) \dd x\right| \\
= \left|\int_{x_2}^{x_1} h(x) \dd x\right| \le |x_1 - x_2|,
$$
so $f$ satisfies a uniform Lipschitz condition of order $1$.  Hence, the class of utility functions that satisfy the Lipschitz condition is at least as large as the class of cumulative distribution functions.  

In general, the bound in Corollary \ref{lipschitz} provides under a specified degree-of-variation on $f$ an upper limit on the increase in the premium $H[\lambda,f(X)]$ resulting from an increase in $\lambda$, the loading parameter.  This enables an insurer to assess the extent to which it is charging suitable additional amounts for perceived increases in risk as measured by higher values of the loading parameter.  
}\end{remark}

\section{Total positivity properties of some actuarial ratios}
\label{actuarialratios}

Concepts of total positivity of higher order are germane to other insurance-related problems.  For $u > 0$ and $v \ge 0$, the well-known {\em upper incomplete gamma function} is defined as 
$$
\Gamma(u,v) = \int_v^\infty x^{u-1} \, \e^{-x} \dd x.
$$
Further, for $c > 0$, define the ratio, 
$$
\mathcal{R}_c(u,v) = \frac{\Gamma(c+u,v)}{\Gamma(u,v)};
$$
this function was shown by \citet{furman_zitikis08b,furman_zitikis08c} to arise in the study of losses from collections of insurable risks, and their Proposition 2.1 proved that $\mathcal{R}_c(u,v)$ is strictly increasing in $u$ for each fixed $v$ and $c$.  Extending this observation, we obtain the following total positivity properties of the function $\mathcal{R}_c$.

\begin{proposition}
\label{Rfunction}
(i) For fixed $v \ge 0$, the function $(c,u) \mapsto \mathcal{R}_c(u,v)$, $c > 0$, $u > 0$, is STP$_\infty$.  

\noindent
(ii) For fixed $u > 0$, the function $(c,v) \mapsto \mathcal{R}_c(u,v)$, $c > 0$, $v \ge 0$, is STP$_\infty$.  

\noindent
(iii)  For fixed $c > 0$, the function $(u,v) \mapsto \mathcal{R}_c(u,v)$, $u > 0$, $v \ge 0$, is SRR$_2$.
\end{proposition}

\Pro
For $r \in \N$, and for $c_1 > \cdots > c_r > 0$ and $u_1 > \cdots > u_r > 0$, consider the $r \times r$ determinant, 
\begin{align*}
\det\big(\Gamma(u_j,v) \, \mathcal{R}_{c_i}(u_j,v)\big) &= \det\big(\Gamma(c_i+u_j,v)\big) \\
&= \det\left( \int_v^\infty x^{c_i+u_j-1} \, \e^{-x} \dd x \right).
\end{align*}
Applying the continuous version of the Binet-Cauchy formula, (\ref{binetcauchy_continuous}), with $\phi_i(x) = x^{c_i}$, $\psi_j(x) = x^{u_j}$, and $\dd\nu(x) = x^{-1} \e^{-x} \dd x$, we find that %the latter determinant equals 
\begin{equation}
\label{Rdeterminant2}
\det\big(\Gamma(u_j,v) \, \mathcal{R}_{c_i}(u_j,v)\big) = \mathop{\idotsint}_{x_1 > \cdots > x_r > v} \, \det\big(x_i^{c_j}\big) \, \det\big(x_i^{u_j}\big) \prod_{j=1}^r x_j^{-1} \e^{-x_j} \dd x_j.
\end{equation}
As shown below in (\ref{powerfunction}), the determinant $\det\big(x_i^{k_j}\big)$ is positive for $x_1 > \cdots > x_r > 0$ and $k_1 > \cdots > k_r \ge 0$.  Therefore the integrand on the right-hand side of (\ref{Rdeterminant2}) is positive on an open subset of $\R^r$, and so the integral is positive.  Hence the determinant on the left-hand side of (\ref{Rdeterminant2}) is positive, and by extracting the factors $\Gamma(u_j,v)$ from that determinant, we obtain 
\begin{equation}
\label{Rdeterminant}
\det\big(\mathcal{R}_{c_i}(u_j,v)\big) > 0.
\end{equation}
Since $r$ was chosen arbitrarily then it follows that the function $(c,u) \mapsto \mathcal{R}_c(u,v)$, $c > 0$, $u > 0$, is STP$_\infty$.  This completes the proof of (i).  

To prove (ii), let ${\bf 1}(x > v)$ denote the indicator function of the interval $(v,\infty)$; then, 
$$
\Gamma(u,v) \mathcal{R}_{c}(u,v) = \int_0^\infty x^c \, {\bf 1}(x > v) \, x^{u-1} e^{-x} \dd x.
$$
For $c_1 > \cdots > c_r > 0$ and $v_1 > \cdots > v_r \ge 0$, it follows from the Binet-Cauchy formula (\ref{binetcauchy_continuous}) that 
\begin{equation}
\label{Rdeterminant3}
\det\big(\Gamma(u,v_j) \, \mathcal{R}_{c_i}(u,v_j)\big) = \mathop{\idotsint}_{x_1 > \cdots > x_r > 0} \, \det\big(x_i^{c_j}\big) \, \det\big({\bf 1}(x_i > v_j)\big) \, \prod_{j=1}^r x_j^{u-1} \e^{-x_j} \dd x_j.
\end{equation}
As noted before, $\det\big(x_i^{c_j}\big) > 0$ for $x_1 > \cdots > x_r > 0$ and $c_1 > \cdots > c_r > 0$.  Further, we note in (\ref{tpindicator}) that for $v_1 > \cdots > v_r \ge 0$, $\det\big({\bf 1}(x_i > v_j)\big) \ge 0$ for all $x_1 > \cdots > x_r > 0$ and also is strictly positive on an open set in $\R^r$.  Therefore, the integrand in (\ref{Rdeterminant3}) is positive on an open set in $\R^r$, so we deduce that the function $(c,v) \mapsto \mathcal{R}_c(u,v)$, $c > 0$, $v \ge 0$, is STP$_\infty$.  

To establish (iii), we apply the ``$2m$-function theorem'' of \citet{rinott_saks}.  For $x > 0$, define the eight functions, 
$$
\begin{array}{ll}
f_1(x) = x^{c+u_1} \, {\bf 1}(x > v_1), \qquad & g_1(x) = x^{c+u_1} \, {\bf 1}(x > v_2) \\[4pt]
f_2(x) = x^{c+u_2} \, {\bf 1}(x > v_2), \qquad & g_2(x) = x^{c+u_2} \, {\bf 1}(x > v_1) \\[4pt]
f_3(x) = x^{u_1} \, {\bf 1}(x > v_2),   \qquad & g_3(x) = x^{u_1} \, {\bf 1}(x > v_1) \\[4pt]
f_4(x) = x^{u_2} \, {\bf 1}(x > v_1),   \qquad & g_4(x) = x^{u_2} \, {\bf 1}(x > v_2)
\end{array}
$$
and, as before, let $\dd\nu(x) = x^{-1} \e^{-x} \dd x$.

We now verify that these functions satisfy the hypotheses of Theorem 1.1 of \cite{rinott_saks}, viz., for $x_1 > x_2 > x_3 > x_4$, 
\begin{equation}
\label{8functionineq}
\prod_{j=1}^4 f_j(x_j) \le \prod_{j=1}^4 g_j(x_j).
\end{equation}
By (\ref{tpindicator}), 
$$
{\bf 1}(z_1 > w_2) \, {\bf 1}(z_2 > w_1) \le {\bf 1}(z_1 > w_1) \, {\bf 1}(z_2 > w_2) 
$$
whenever $z_1 > z_2$ and $w_1 > w_2$.  Applying this result repeatedly, and noting that $v_1 > v_2$, we obtain 
\begin{align*}
{\bf 1}(x_1 > v_1) \, {\bf 1}(x_2 > v_2) \, {\bf 1}&(x_3 > v_2) \, {\bf 1}(x_4 > v_1) \\
&\le {\bf 1}(x_1 > v_1) \, {\bf 1}(x_2 > v_2) \, {\bf 1}(x_3 > v_1) \, {\bf 1}(x_4 > v_2) \\
& \le {\bf 1}(x_1 > v_1) \, {\bf 1}(x_2 > v_1) \, {\bf 1}(x_3 > v_2) \, {\bf 1}(x_4 > v_2) \\
& \le {\bf 1}(x_1 > v_2) \, {\bf 1}(x_2 > v_1) \, {\bf 1}(x_3 > v_1) \, {\bf 1}(x_4 > v_2).
\end{align*}
Multiplying each side of this inequality by $x_1^{c+u_1} \, x_2^{c+u_2} \, x_3^{u_1} \, x_4^{u_2}$, we obtain (\ref{8functionineq}); moreover, that inequality is strict on an open set in $(0,\infty)^4$.  It is also trivial that the measure $\dd\nu(x)$ is an {\it FKG measure} (\citet[p. 270]{rinott_saks}), so we obtain 
\begin{align*}
\Gamma(c+u_1,v_1) \Gamma(c+u_2,v_2) \Gamma(u_1,v_2) \Gamma(u_2,v_1) &= \prod_{j=1}^4 \int_0^\infty f_j(x) \dd\nu(x) \\
&< \prod_{j=1}^4 \int_0^\infty g_j(x) \dd\nu(x) \\
&= \Gamma(c+u_1,v_2) \Gamma(c+u_2,v_1) \Gamma(u_1,v_1) \Gamma(u_2,v_2).
\end{align*}
Dividing both sides of this inequality by $\Gamma(u_1,v_1) \Gamma(u_2,v_2) \Gamma(u_1,v_2) \Gamma(u_2,v_1)$, we obtain $\det\big(\mathcal{R}_c(u_i,v_j)\big) < 0$ for $u_1 > u_2$ and $v_1 > v_2$.  Hence, $\mathcal{R}_c(u,v)$ is SRR$_2$ in $(u,v)$.  
$\qed$

\begin{remark}{\rm
(i) As a special case of (\ref{Rdeterminant}), suppose that $r = 2$, $c_1 = c$, and $c_2 = 0$; since $\mathcal{R}_0(u,v) = 1$ then (\ref{Rdeterminant}) reduces to the monotonicity result of \citet{furman_zitikis08b}.  More generally, Proposition \ref{Rfunction} can be applied to obtain inequalities for the higher moments of sums of risks similar to the way in which higher moment inequalities are described in Remark \ref{remarkhighermoments}.  

(ii) We remark that the function $(u,v) \mapsto \mathcal{R}_c(u,v)$, $u > 0$, $v \ge 0$ is not RR$_3$.  For $c = 3.5$, $(u_1,u_2,u_3) = (4,3,2)$, and $(v_1,v_2,v_3) = (6,5,4)$, we calculate that $\det\big(\mathcal{R}_c(u_i,v_j)\big) = 7.04$ 
%7.03758272
which, since it is positive, violates the RR$_3$ condition.  Here and throughout, all numerical computations were carried out in high precision, with accuracy to over $100$ significant digits.  
}\end{remark}

\medskip

We will also provide some results on a more complex ratio defined by \citet{furman_zitikis08c}.  Define, for $u > 0$ and $v \ge 0$ the function
\begin{equation}
\label{Qfunction}
Q(u,v) = \frac{\Gamma(u,v)}{\Gamma(u)} = \int_v^\infty \, \frac{1}{\Gamma(u)} \, x^{u-1} \e^{-x} \dd x.
\end{equation}
For fixed $u > 0$, the function $v \mapsto Q(u,v)$, $v \ge 0$, is strictly decreasing, so we denote by $v \mapsto Q^{-1}(u,v)$, $0 \le v \le 1$, the corresponding inverse function.  \citet{furman_zitikis08c} defined the function, 
$$
\mathcal{C}_c(u,v) = Q\big(c+u,Q^{-1}(u,v)\big) = \frac{\Gamma\big(c+u,Q^{-1}(u,v)\big)}{\Gamma(c+u)},
$$
$c > 0$, $u > 0$, $0 \le v \le 1$; further, they proved that, for fixed $c$ and $v$, the function $u \mapsto \mathcal{C}_c(u,v)$ is decreasing.  In light of Proposition \ref{Rfunction} we will investigate the total positivity properties of $\mathcal{C}_c(u,v)$, obtaining the following result.  

\begin{proposition}
\label{Cfunction}
(i) For fixed $v \in [0,1]$, the function $(c,u) \mapsto \Gamma(c+u) \, \mathcal{C}_c(u,v)$, $c > 0$, $u > 0$, is STP$_2$.  

\noindent
(ii) For fixed $u > 0$, the function $(c,v) \mapsto \mathcal{C}_c(u,v)$, $c > 0$, $0 \le v \le 1$, is SRR$_\infty$.
\end{proposition}

\Pro 
We have 
\begin{align*}
\Gamma(c+u) \, \mathcal{C}_c(u,v) &= \Gamma\big(c+u,Q^{-1}(u,v)\big) \\
&= \int_{Q^{-1}(u,v)}^\infty x^{c+u-1} \, e^{-x} \dd x \\
&= \int_0^\infty x^{c+u-1} \, {\bf 1}\big(x > Q^{-1}(u,v)\big) \, e^{-x} \dd x.
\end{align*}
For $c_1 > c_2 > 0$ and $u_1 > u_2 > 0$, it follows from the Binet-Cauchy formula that 
\begin{multline}
\label{bcfgammacu}
\det\big(\Gamma(c_i+u_j) \, \mathcal{C}_{c_i+u_j}(u_j,v)\big) \\
= \mathop{\idotsint}_{x_1 > x_2 > 0} \, \det(x_i^{c_j}) \, \det\Big(x_i^{u_j} {\bf 1}\big(x_i > Q^{-1}(u_j,v)\big)\Big) \prod_{j=1}^2 x_j^{-1} e^{-x} \dd x_j.
\end{multline}
For $2 \times 2$ matrices $A = (a_{ij})$ and $B = (b_{ij})$, there holds the identity, 
\begin{equation}
\label{2by2detidentity}
\det(a_{ij} b_{ij}) = b_{11} b_{22} \det(A) + a_{12} a_{21} \det(B).
\end{equation}
We set $a_{ij} = x_i^{u_j}$; then $\det(A) = \det(x_i^{c_j}) > 0$ and $a_{12} a_{21} > 0$ for $x_1 > x_2 > 0$ and $c_1 > c_2 > 0$.  Also, let $b_{ij} = {\bf 1}\big(x_i > Q^{-1}(u_j,v)\big)$; then we need to determine the sign of 
\begin{equation}
\label{detB}
\det(B) = \det\Big({\bf 1}\big(x_i > Q^{-1}(u_j,v)\big)\Big)
\end{equation}
for $x_1 > x_2 > 0$ and $u_1 > u_2 > 0$.  

We claim that if $u_1 > u_2 > 0$ then $Q^{-1}(u_1,v) > Q^{-1}(u_2,v)$ for all $v > 0$.  By applying $Q(u_1,\cdot)$ to both sides of this inequality, and noting that $v \mapsto Q(u_1,v)$ is decreasing, we see that the claim is equivalent to $v < Q\big(u_1,Q^{-1}(u_2,v)\big)$; and by replacing $v$ further by $Q(u_2,v)$, the claim is now seen to be equivalent to $Q(u_2,v) < Q(u_1,v)$, $v > 0$.  So, consider 
\begin{align*}
\Gamma(u_1) \Gamma(u_2) \big[Q(u_1,v) - Q(u_2,v)\big] &= \Gamma(u_2) \Gamma(u_1,v) - \Gamma(u_1) \Gamma(u_2,v) \\
&= \int_{x_2=0}^\infty \int_{x_1=v}^\infty (x_1^{u_1} x_2^{u_2} - x_1^{u_2} x_2^{u_1}) \prod_{j=1}^2 x_j^{-1} e^{-x_j} \dd x_j \\
&= \int_{x_2=0}^\infty \int_{x_1=v}^\infty \det(x_i^{u_j}) \prod_{j=1}^2 x_j^{-1} e^{-x_j} \dd x_j.
\end{align*}
Decomposing the interval $(0,\infty)$ into $(0,v) \cup [v,\infty)$, we obtain 
\begin{multline*}
\Gamma(u_1) \Gamma(u_2) \big[Q(u_1,v) - Q(u_2,v)\big] \\
= \int_v^\infty \int_v^\infty \det(x_i^{u_j}) \prod_{j=1}^2 x_j^{-1} e^{-x_j} \dd x_j + \int_0^v \int_v^\infty \det(x_i^{u_j}) \prod_{j=1}^2 x_j^{-1} e^{-x_j} \dd x_j.
\end{multline*}
By anti-symmetry, 
$$
\int_v^\infty \int_v^\infty \det(x_i^{u_j}) \prod_{j=1}^2 x_j^{-1} e^{-x_j} \dd x_j = 0,
$$
therefore
\begin{align*}
\Gamma(u_1) \Gamma(u_2) \big[Q(u_1,v) - Q(u_2,v)\big] &= \int_0^v \int_v^\infty \det(x_i^{u_j}) \prod_{j=1}^2 x_j^{-1} e^{-x_j} \dd x_j \\
&\equiv \mathop{\int\int}_{x_1 > v > x_2 > 0} \det(x_i^{u_j}) \prod_{j=1}^2 x_j^{-1} e^{-x_j} \dd x_j.
\end{align*}
Since $\det(x_i^{u_j}) > 0$ for $x_1 > x_2 > 0$ and $u_1 > u_2 > 0$ then we obtain $Q(u_1,v) > Q(u_2,v)$, equivalently, $Q^{-1}(u_1,v) > Q^{-1}(u_2,v)$.

Since $u_1 > u_2$ implies $Q^{-1}(u_1,v) > Q^{-1}(u_2,v)$ then it follows from (\ref{detB}) and (\ref{tpindicator}) that for $x_1 > x_2$ and $u_1 > u_2$, 
$$
\det(B) %\equiv \det\Big({\bf 1}\big(x_i > Q^{-1}(u_j,v)\big)\Big) 
= \begin{cases}
1, & x_1 > Q^{-1}(u_1,v) > x_2 > Q^{-1}(u_2,v) \\
0, & \hbox{otherwise}
\end{cases} \ .
$$
Therefore for $u_1 > u_2 > 0$, $\det(B)$, $b_{11}$, and $b_{22}$, as functions of $(x_1,x_2)$, are nonnegative everywhere and are positive on an open subset of the orthant $\{x_1 > x_2 > 0\}$.  Applying (\ref{2by2detidentity}), we find that 
$$
\det\Big(x_i^{u_j} {\bf 1}\big(x_i > Q^{-1}(u_j,v)\big) \Big) > 0
$$
on an open subset of the orthant $\{x_1 > x_2 > 0\}$, consequently, the same applies to the integrand in (\ref{bcfgammacu}), so the integral is positive.  Therefore, the function $(c,u) \mapsto \Gamma\big(c+u,Q^{-1}(u,v)\big)$, $c > 0$, $u > 0$, is STP$_2$.  The proof of (i) now is complete.  

To prove (ii), let $r \in \N$, $c_1 > \cdots > c_r > 0$, and $v_1 > \cdots > v_r \ge 0$; then, 
\begin{align*}
\Big[\prod_{i=1}^r \Gamma(c_i+u)\Big] \det\big(\mathcal{C}_{c_i}(u,v_j)\big) &= \det\big(\Gamma(c_i+u) \, \mathcal{C}_{c_i}(u,v_j)\big) \\
&= \det\Big(\Gamma\big(c_i+u,Q^{-1}(u,v_j)\big)\Big) \\
&= \det\Big(\int_{Q^{-1}(u,v_j)}^\infty x^{c_i+u-1} \, e^{-x} \, \dd x\Big).
\end{align*}
Applying the Binet-Cauchy formula, we obtain 
\begin{multline}
\label{cvtoCcuv}
\Big[\prod_{i=1}^r \Gamma(c_i+u)\Big] \det\big(\mathcal{C}_{c_i}(u,v_j)\big) \\
= \mathop{\idotsint}_{x_1 > \cdots > x_r > 0} \, \det(x_i^{c_j}) \, \det\Big({\bf 1}\big(x_i > Q^{-1}(u,v_j)\big) \Big) \prod_{j=1}^r x_j^{u-1} e^{-x} \dd x_j.
\end{multline}
As before, $\det(x_i^{c_j}) > 0$ for $x_1 > \cdots > x_r > 0$.  To derive the sign of the remaining determinant, note that $Q^{-1}(u,v_r) > Q^{-1}(u,v_{r-1}) > \cdots > Q^{-1}(u,v_r)$ for $v_1 > \cdots > v_r$.  Therefore, for $x_1 > \cdots > x_r$ and $v_1 > \cdots > v_r$, 
\begin{align*}
\det&\Big({\bf 1}\big(x_i > Q^{-1}(u,v_j)\big)\Big) \\
&= (-1)^{r(r-1)/2} \det\Big({\bf 1}\big(x_i > Q^{-1}(u,v_{r-j+1})\big)\Big) \\
&= \begin{cases}
(-1)^{r(r-1)/2}, & x_1 > Q^{-1}(u,v_r) > x_2 > Q^{-1}(u,v_{r-1}) > \cdots > x_r > Q^{-1}(u,v_1) \\
0, & \hbox{otherwise}
\end{cases}\ .
\end{align*}
It now follows that the sign of (\ref{cvtoCcuv}) is $(-1)^{r(r-1)/2}$.  Since $r$ was chosen arbitrarily, it follows that the function $\mathcal{C}_c(u,v)$ is SRR$_\infty$ in $(c,v)$.  
$\qed$

\begin{remark}{\rm 
As we noted before, \citet{furman_zitikis08c} proved that the function $u \mapsto \mathcal{C}_c(u,v)$ is decreasing, a result which raises the issue of whether the function $(c,u) \mapsto \mathcal{C}_c(u,v)$, $c > 0$, $u > 0$ is RR$_2$.  

Let $F_u$ denote the cumulative distribution function of a gamma-distributed random variable having index parameter $u$, i.e., a random variable whose probability density function is the integrand in (\ref{Qfunction}).   Then, $\mathcal{C}_c(u,v) = 1 - F_{c+u}\big(Q^{-1}(u,v)\big)$.  Since $Q(u,v) = 1- F_u(v)$ then $Q^{-1}(u,v) = F_u^{-1}(1-v)$, so we obtain 
$$
\mathcal{C}_c(u,v) = 1 - F_{c+u}\big(F_u^{-1}(1-v)\big).
$$
We performed extensive calculations using this identity and determined that the function $(c,u) \mapsto \mathcal{C}_c(u,v)$ is not RR$_2$ as many of its $2 \times 2$ determinants are positive.  Further, this function is not TP$_2$, although it seems to fail barely to be so; for $v = 0.211$, $(c_1,c_2) = (4.047,1.210)$, and $(u_1,u_2) = (3.203,0.189)$, we obtained $\det\big(\mathcal{C}_{c_i}(u_j,v)\big) = -0.026$, and this negative value of the determinant was similar in magnitude to all other negative values that we found.  

We also carried out calculations regarding the total positivity properties of the function $(u,v) \mapsto \mathcal{C}_c(u,v)$ and found substantial evidence that this function is both STP$_2$ and STP$_3$.  As we have not been able to establish such results analytically, we pose them as open problems.
}\end{remark}

\section{Seven classes of weight functions}
\label{sevenweightfunctions}

In this section, we determine the total positivity properties of some classes of weight functions treated by \citet[Section 3]{sendov_etal}.  

\begin{example}
\label{example:esscher}
{\rm 
Let $w_1(\lambda,x) = \e^{\lambda x}$.  The corresponding weighted premium, $H[\lambda,X]$, is called the {\it Esscher premium}; see \citet{sendov_etal} and references given therein.  It is well-known that the weight function $w_1$ is STP$_\infty$ \citep[p.~15]{karlin}.  

Indeed, by a result of \citet[p. 233]{gross_richards}, for each $r \ge 2$, the $r \times r$ determinant, $\det\big(w_1(\lambda_i,x_j)\big)$ has an integral representation, 
$$
\frac{\det\big(w_1(\lambda_i,x_j)\big)}{\operatornamewithlimits\prod\limits_{1 \le i < j \le r} (\lambda_i - \lambda_j)(x_i - x_j)} = \int_{U} \Phi(\Lambda,X,u) \dd\nu(u),
$$
where $U$ is a certain set of $r \times r$ matrices, $\Lambda = (\lambda_1,\ldots,\lambda_r)$, $X = (x_1,\ldots,x_r)$, $\Phi(\Lambda,X,u)$ is a strictly positive function, and $\nu$ is a probability measure on $U$.  This integral formula yields immediately the positivity of the determinant $\det\big(w_1(\lambda_i,x_j)\big)$ for $\lambda_1 > \cdots > \lambda_r$ and $x_1 > \cdots > x_r$.  Hence $w_1$ is STP$_r$ for all $r$ and therefore also STP$_\infty$.  

More generally, if $F:\R \to \R$ is strictly increasing then the weight function $\tilde{w}_1(\lambda,x) = \exp\big(\lambda F(x)\big)$ is STP$_\infty$, and the corresponding premium is known as the {\em Aumann-Shapley premium} \citep{furman_zitikis09}.  With regard to the closing comments in Remark \ref{remarkhighermoments}, we would advise an insurer to base its premium calculations on an STP$_\infty$ weight function if the possible values of the loss variable $X$ are greatly dispersed, i.e., if $X$ has an extremely large variance-to-mean ratio.
}\end{example}

A consequence of the STP$_\infty$ property of $w_1$ is that the weight function $w_1(\log \lambda,x) = \lambda^x$, $(\lambda,x) \in \R_+ \times \R$, is STP$_\infty$.  That is, for any $r \in \N$, the $r \times r$ determinant, 
\begin{equation}
\label{powerfunction}
\det\big(\lambda_i^{x_j}\big) > 0
\end{equation}
for $\lambda_1 > \cdots > \lambda_r > 0$ and $x_1 > \cdots > x_r$.  This result holds because the transformation $\lambda \to \log \lambda$ is strictly increasing and therefore preserves the total positivity properties of $w_1(\lambda,x)$.  In the actuarial literature, the weight function $(\lambda,x) \mapsto \lambda^x$, $\lambda, x > 0$, gives rise to a weighted premium known as the {\em size-biased premium} \citep{furman_zitikis09}.  
%We shall apply (\ref{powerfunction}) repeatedly in the sequel.  

\begin{example}
\label{example:ctp}
{\rm 
Let $w_2(\lambda,x) = {\bf 1}(x > \lambda)$, the weight function corresponding to  the {\it conditional tail expectation} (CTE) premium, $H[\lambda,X] = \E(X | X > \lambda)$.  It is well-known that the weight function $w_2$ is TP$_\infty$ \citep[p.~16]{karlin}.  Indeed, for any $r = 1, 2, \ldots$, and for $\lambda_1 > \cdots > \lambda_r$ and $x_1 > \cdots > x_r$, 
\begin{equation}
\label{tpindicator}
\det\big(w_2(\lambda_i,x_j)\big) = 
\begin{cases}
1, & \hbox{if } x_1 > \lambda_1 > x_2 > \lambda_2 > \cdots > x_r > \lambda_r \\
0, & \hbox{otherwise}
\end{cases}
\end{equation}
which proves that the determinant is nonnegative.  
}\end{example}

\begin{example}
\label{example:kamps}
{\rm 
Let $w_3(\lambda,x) = 1 - \e^{-x/\lambda}$.  The corresponding weighted premium $H[\lambda,X]$ is called the {\it Kamps premium} (see \citet{sendov_etal} and the references therein).  The weight function $w_3$ is STP$_\infty$, as we now prove.  

It is straightforward to verify that 
\begin{eqnarray*}
w_3(\lambda,x) & = & \int_0^x \lambda^{-1} \e^{-t/\lambda} \dd t \\
& \equiv & \int_0^\infty w(\lambda,t) \, w_2(t,x) \dd t,
\end{eqnarray*}
where $w(\lambda,t) = \lambda^{-1} \e^{-t/\lambda}$, and $w_2(t,x) = {\bf 1}(x > t)$ is the weight function given in Example \ref{example:ctp}.  Applying the Binet-Cauchy formula (\ref{binetcauchy_continuous}), we obtain for each $r \ge 2$, 
\begin{equation}
\label{bcf-w3}
\det\big(w_3(\lambda_i,x_j)\big) = \idotsint\limits_{t_1 > \cdots > t_r} \det\big(w(\lambda_i,t_j)\big) \, \det\big(w_2(t_i,x_j)\big) \dd t_1 \cdots \dd t_r.
\end{equation}

Suppose that $\lambda_1 > \cdots > \lambda_r$ and $x_1 > \cdots > x_r$.  By Example \ref{example:esscher}, the function $w(\lambda,t)$ is STP$_\infty$, so $\det\big(w(t_i,x_j)\big)$ is positive on the orthant $\{(t_1,\ldots,t_r): t_1 > \cdots > t_r\}$.  Also, by (\ref{tpindicator}), the determinant $\det\big(w_2(t_i,x_j)\big)$ is positive on an open neighborhood in the same orthant.  Therefore, the integrand in (\ref{bcf-w3}) is positive on an open set, hence the integral is positive for any choice of $(\lambda_1,\ldots,\lambda_r)$ and $(x_1,\ldots,x_r)$.  Therefore, $w_3(\lambda,x)$ is STP$_r$ for all $r$, hence it is STP$_\infty$.  

These results for $w_3(\lambda,x)$ also extend to a more general class of weight functions.  For each nonnegative integer $k$, define the weight function 
$$
w_{3,k}(\lambda,x) = k!\bigg[1 - \e^{-x/\lambda} \sum_{j=0}^k \frac{(x/\lambda)^j}{j!}\bigg].
$$
For $k = 0$, $w_{3,k}(\lambda,x)$ reduces to $w_3(\lambda,x)$, the weight function corresponding to Kamps' premium.  By repeated integration-by-parts, we obtain 
\begin{eqnarray*}
w_{3,k}(\lambda,x) & = & \lambda^{-(k+1)} \int_0^x t^k \e^{-t/\lambda} \dd t \\
& = & \int_0^\infty w(\lambda,t) \, w_2(t,x) \dd t,
\end{eqnarray*}
where $w(\lambda,t) = \lambda^{-(k+1)} \e^{-t/\lambda}$ and $w_2(t,x)$ is the weight function in Example \ref{example:ctp}.  
Finally, we proceed using arguments similar to the case of $w_3$: Since $w(\lambda,t)$ and $w_2(\lambda,t)$ are TP$_\infty$ then, by applying the Basic Composition Formula (\ref{bcf_continuous}) and the Binet-Cauchy formula (\ref{binetcauchy_continuous}), we deduce that $w_{3,k}(\lambda,x)$ is STP$_\infty$.  
}\end{example}

\begin{example}
\label{example:pseudopoisson}
{\rm 
The fourth weight function considered by \citet[Section 3]{sendov_etal} is
%\begin{equation}
%\label{w4lambdax}
$$
\wtilde{w}_4(\lambda,x) = \exp\bigg(\frac{(1+x)^\lambda - 1}{\lambda}\bigg) - x,
$$
%\end{equation}
$\lambda, x > 0$.  We replace $x$ by $\e^x - 1$, a transformation that is strictly increasing and therefore preserves any total positivity properties of $\tilde{w}_4(\lambda,x)$.  Then we are to determine the total positivity properties of 
$$
w_4(\lambda,x) := \wtilde{w}_4(\lambda,\e^x - 1) = \exp\big(f(\lambda,x)\big) - \e^x + 1,
$$
$\lambda, x > 0$, where 
$$
f(\lambda,x) = \frac{\e^{\lambda x} - 1}{\lambda}.
$$
We also define $f(0,x)$ by right-continuity: 
$$
f(0,x) := \lim_{\lambda \to 0+} f(\lambda,x) = x.
$$
Then, 
\begin{equation}
\label{w4function}
w_4(\lambda,x) = \exp\big(f(\lambda,x)\big) - \exp\big(f(0,x)\big) + 1.
\end{equation}
For $r \in \N$, let $\lambda_1 > \cdots > \lambda_r > 0$, $x_1 > \cdots > x_r > 0$, and consider the $r \times r$ determinant, 
$$
\det\big(w_4(\lambda_i,x_j)\big) = \left|\begin{array}{ccc}
w_4(\lambda_1,x_1) & \cdots & w_4(\lambda_1,x_r) \\
w_4(\lambda_2,x_1) & \cdots & w_4(\lambda_2,x_r) \\
\vdots & \vdots & \vdots \\
w_4(\lambda_{r-1},x_1) & \cdots & w_4(\lambda_{r-1},x_r) \\
w_4(\lambda_r,x_1) & \cdots & w_4(\lambda_r,x_r)
\end{array}\right|.
$$
For $i=1,\ldots,r-1$, we subtract row $i+1$ from row $i$, obtaining 
\begin{align*}
\det\big(w_4(\lambda_i,x_j)\big) = D_1 + D_2,
\end{align*}
where 
$$
D_1 = \left|\begin{array}{ccc}
w_4(\lambda_1,x_1) - w_4(\lambda_2,x_1) & \cdots & w_4(\lambda_1,x_r) - w_4(\lambda_2,x_r) \\
w_4(\lambda_2,x_1) - w_4(\lambda_3,x_1) & \cdots & w_4(\lambda_2,x_r) - w_4(\lambda_3,x_1)\\
\vdots & \vdots & \vdots \\
w_4(\lambda_{r-1},x_1) - w_4(\lambda_r,x_1)& \cdots & w_4(\lambda_{r-1},x_r) - w_4(\lambda_r,x_1) \\
w_4(\lambda_r,x_1) - w_4(0,x_1) & \cdots & w_4(\lambda_r,x_r) - w_4(0,x_1)
\end{array}\right| \\
$$
and 
\begin{equation}
\label{D2}
D_2 = 
\left|\begin{array}{ccc}
w_4(\lambda_1,x_1) - w_4(\lambda_2,x_1) & \cdots & w_4(\lambda_1,x_r) - w_4(\lambda_2,x_r) \\
w_4(\lambda_2,x_1) - w_4(\lambda_3,x_1) & \cdots & w_4(\lambda_2,x_r) - w_4(\lambda_3,x_1)\\
\vdots & \vdots & \vdots \\
w_4(\lambda_{r-1},x_1) - w_4(\lambda_r,x_1)& \cdots & w_4(\lambda_{r-1},x_r) - w_4(\lambda_r,x_1) \\
w_4(0,x_1) & \cdots & w_4(0,x_1)
\end{array}\right|.
\end{equation}

Define 
$$
w_{41}(\lambda,x) = \frac{\partial}{\partial \lambda} \omega_4(\lambda,x) \equiv \frac{\partial}{\partial \lambda} \exp\big(f(\lambda,x)\big),
$$
and set $\lambda_{r+1} \equiv 0$.  By Taylor's theorem, there exists $\rho_i \in (\lambda_{i+1},\lambda_i)$ such that 
\begin{equation}
\label{taylor}
w_4(\lambda_i,x) - w_4(\lambda_{i+1},x) = (\lambda_i - \lambda_{i+1}) w_{41}(\rho_i,x),
\end{equation}
$i=1,\ldots,r$.  Therefore, 
\begin{align*}
D_1 &= \det\big((\lambda_i - \lambda_{i+1}) w_{41}(\rho_i,x_j)\big) \\
&\equiv \prod_{i=1}^r (\lambda_i - \lambda_{i+1}) \cdot \det\big(w_{41}(\rho_i,x_j)\big),
\end{align*}
where $\lambda_1 > \rho_1 > \lambda_2 > \rho_2 > \cdots > \lambda_r > \rho_r > 0$.  So, to prove that $D_1$ is positive, it suffices to show that $w_{41}$ is STP$_r$, and we begin by observing from (\ref{w4function}) that 
\begin{align}
\label{w41asderiv}
\det\big(w_{41}(\lambda_i,x_j)\big) &= \det\Big(\frac{\partial}{\partial \lambda_i} \exp\big(f(\lambda_i,x_j)\big)\Big) \nonumber \\
&= \frac{\partial^r}{\partial \lambda_1 \cdots \partial \lambda_r} \det\Big(\exp\big(f(\lambda_i,x_j)\big)\Big).
\end{align}

We now recall the {\it Bell} (or {\it exponential}\,) {\it polynomials} $B_k$, $k=0,1,2,\ldots$, defined through the generating function, 
\begin{equation}
\label{Bellgeneratingfunction}
\exp\big(u(\e^t - 1)\big) = \sum_{k=0}^\infty B_k(u) \frac{t^k}{k!}.
\end{equation}
We refer to \citet[p.~133 ff.]{comtet} and \citet[pp.~63-67]{roman} for further details on these polynomials.  For $k \ge 1$, $B_k(u)$ is monic and of degree $k$; moreover, 
\begin{equation}
\label{BellStirling}
B_k(u) = \sum_{m=1}^k S(k,m) u^m,
\end{equation}
where the coefficients $S(k,m)$ are the {\it Stirling numbers of the second kind}, viz., the number of partitions of a set of size $k$ into $m$ non-empty subsets \citep[p.~50]{comtet}.  In particular, $S(k,1) = S(k,k) = 1$, and $S(k,m) = 0$ if $m > k$.  
%The simplest examples of the Bell polynomials are: $B_0(u) = 1$, $B_1(u) = u$, and $B_2(u) = u^2 + u$.  

An alternative representation for the Bell polynomials arises from the observation that the left-hand side of (\ref{Bellgeneratingfunction}) is, for $u > 0$, the moment-generating function of $U$, a Poisson-distributed random variable with mean parameter $u$; therefore, 
\begin{equation}
\label{poissonmoments}
B_k(u) = E(U^k) = \sum_{m=0}^\infty \frac{\e^{-u} u^m}{m!} m^k.
\end{equation}
We now apply to (\ref{poissonmoments}) the discrete Binet-Cauchy formula (\ref{binetcauchy_discrete}) with $\phi_i(m) = \e^{-u_i} u_i^m$ and $\psi_k(m) = m^k$, each of which is STP$_\infty$ by (\ref{powerfunction}), and weights $\nu(m) = 1/m!$.  Written explicitly, we have, for $k_1 > \cdots > k_r \ge 0$ and $u_1 > \cdots> u_r > 0$, 
$$
\det\big(B_{k_i}(u_j)\big) = \e^{-(u_1+\cdots+u_r)} \sum_{m_1 > \cdots > m_r \ge 0} \frac{1}{m_1! \cdots m_r!} \det(u_i^{m_j}) \det(m_i^{k_j}).
$$
Then the positivity of $\det\big(B_{k_i}(u_j)\big)$ follows from the positivity of each determinant inside the summation.  Since $r$ was chosen arbitrarily then it follows that $B_k(u)$ is STP$_\infty$ in $(k,u)$.  

Define 
\begin{equation}
\label{Pkpolynomial}
\wtilde{B}_k(\lambda) := \lambda^k B_k(\lambda^{-1}) = \sum_{m=0}^{k-1} S(k,k-m) \lambda^m.
\end{equation}
%Since $S(k,1) = 1$ then the polynomial $\wtilde{B}_k$ is monic and of degree $k-1$.
Then, by (\ref{Bellgeneratingfunction}) and (\ref{BellStirling}), 
\begin{align}
\label{Pkgeneratingfunction}
\exp\big(f(\lambda,x)\big) = \sum_{k=0}^\infty B_k(\lambda^{-1}) \frac{(\lambda x)^k}{k!} = \sum_{k=0}^\infty \wtilde{B}_k(\lambda) \frac{x^k}{k!}.
\end{align}
Applying to (\ref{Pkgeneratingfunction}) the Binet-Cauchy formula (\ref{binetcauchy_continuous}), we obtain for $\lambda_1>\cdots>\lambda_r>0$ and $x_1>\cdots>x_r>0$, 
\begin{align}
\label{detw41D1}
\det\Big(\exp\big(f(\lambda_i,x_j)\big)\Big) &= \det\bigg(\sum_{k=0}^\infty \frac{1}{k!} \, \wtilde{B}_k(\lambda_i) \, x_j^k\bigg) \nonumber \\
&= \sum_{k_1 > \cdots > k_r \ge 0} \frac{1}{k_1! \cdots k_r!} \, \det\big(\wtilde{B}_{k_j}(\lambda_i)\big) \, \det\big(x_i^{k_j}\big).
\end{align}
By (\ref{powerfunction}), $\det\big(x_i^{k_j}\big) > 0$ for $x_1 > \cdots > x_r$ and $k_1 > \cdots > k_r$.  

We note two consequences of (\ref{Pkpolynomial}).  First, since $S(k,1) = 1$ then the polynomial $\wtilde{B}_k$ is monic and of degree $k-1$.  Second, since $S(k,m) = 0$ for $m > k$ and $S(k,k) = 1$ then the polynomials $\{\wtilde{B}_k(\lambda): k=0,1,2,\ldots\}$ satisfy a linear system of equations in terms of $\{\lambda^m: m=0,1,2,\ldots\}$, with a triangular matrix of coefficients having the $(k,m)$th entry equal to $S(k,m)$, $1 \le m \le k$.  Writing out these equations for rows $k_r,k_{r-1},\ldots,k_1$, in that order, and for columns $1,2,\ldots,r$  
%$\lambda_1 > \cdots > \lambda_r$ and $k_r < k_{r-1} < \cdots < k_1$ 
results in a matrix equation, 
\begin{equation}
\label{BS}
\wtilde{\bB} = \bS \bLambda,
\end{equation}
where the $r \times r$ matrix $\wtilde{\bB}$ has $(i,j)$th entry $\wtilde{B}_{k_{r-i+1}}(\lambda_j)$, $i,j=1,\ldots,r$; $\bS$ is $r \times k_1$ with $(i,j)$th entry $S(k_{r-i+1},k_{r-i+1}-j+1)$, $i=1,\ldots,r$, $j=1,\ldots,k_1$; and $\bLambda$ is $k_1 \times r$ with $(i,j)$th entry $\lambda_j^{i-1}$, $i=1,\ldots,k_1$, $j=1,\ldots,r$.  

Each $r \times r$ minor of $\bLambda$, being of the form $\det\big(\lambda_j^{l_i}\big)$ with $\lambda_1 > \cdots > \lambda_r$ and $l_1 < \cdots < l_r$, is non-zero and has sign $(-1)^{r(r-1)/2}$ as $r(r-1)/2$ row interchanges are needed to order the $\lambda_i$ and $l_i$ similarly.  

To determine the sign of the minors of $\bS$, we apply the results of \cite[Section~5]{brenti}; cf.~\cite{mongelli}.  According to those results, the infinite matrix $S(k,m)$ is totally positive, i.e., all minors of the matrix $S(k,m)$, where $k$ and $m$ are similarly ordered, are nonnegative.  In the $i$th row of $\bS$, the columns are indexed by the decreasing sequence $k_{r-i+1}-j+1$, $j=1,\ldots,k_1$; and in the $j$th column, the rows are indexed by the increasing sequence $k_{r-i+1}$, $i=1,\ldots,r$; therefore, it follows that each non-zero $r \times r$ minor of $\bS$ also has sign $(-1)^{r(r-1)/2}$.  Further, if $k_1,\ldots,k_r$ are consecutive integers then the resulting matrix is lower triangular with non-zero diagonal entries, so the corresponding minor of $\bS$ is non-zero.  

By the classical Binet-Cauchy formula, $\det(\wtilde{\bB})$ equals a sum of products of $r \times r$ minors of $\bS$ and $\bLambda$ \cite[p.~1]{karlin}.  By the preceding discussion, each such product is nonnegative; this establishes the positivity of each $r \times r$ minor of $\wtilde{\bB}$, proving that it is at least TP$_r$.  Further, the sum of all such products of minors is positive since some minors of $\bS$, and all minors of $\bLambda$, are non-zero.  Therefore, $\wtilde{\bB}$ is STP$_r$; and since $r$ is arbitrary then it is STP$_\infty$.  

To complete the proof that $w_{41}(\lambda,x)$ is STP$_r$, we need to show that (\ref{w41asderiv}) is positive.  By the same argument as at (\ref{BS}), {\it infra}, we find that 
\begin{equation}
\label{Btildelambdaderivs}
\frac{\partial^r}{\partial \lambda_1 \cdots \partial \lambda_r} \det\big(\wtilde{B}_{k_j}(\lambda_i)\big)
\end{equation}
is a sum of products of minors of $\bS$ with derivatives of minors of $\bLambda$.  However, the derivatives of the minors of $\bLambda$ are the form 
\begin{equation}
\label{vandermondederivs}
\frac{\partial^r}{\partial \lambda_1 \cdots \partial \lambda_r} \det\big(\lambda_j^{l_i}\big) \equiv \det\Big(\frac{\partial}{\partial \lambda_i} \lambda_j^{l_i}\Big) = \frac{l_1 \cdots l_r}{\lambda_1 \cdots \lambda_r} \det\big(\lambda_j^{l_i}\big),
\end{equation}
which is of the same sign, viz., $(-1)^{r(r-1)/2}$, as each $r \times r$ minor of $\bLambda$.  Therefore, the derivatives (\ref{Btildelambdaderivs}) are nonnegative, and some are positive.  By differentiating the series (\ref{detw41D1}), it follows that $w_{41}(\lambda,x)$ is STP$_r$, and hence is STP$_\infty$.  

For future reference, we note that in addition to (\ref{w41asderiv}) being positive for $\lambda_1 > \cdots > \lambda_r$ and $x_1 > \cdots > x_r$, there also holds 
\begin{equation}
\label{mixedderivdet}
\det\Big(\frac{\partial^2}{\partial x_j \partial \lambda_i} \exp\big(f(\lambda_i,x_j)\big)\Big) > 0
\end{equation}
under the same conditions.  To prove this, we observe that the above determinant equals 
$$
\frac{\partial^r}{\partial x_1 \cdots \partial x_r} \frac{\partial^r}{\partial \lambda_1 \cdots \partial \lambda_r} \det\Big(\exp\big(f(\lambda_i,x_j)\big)\Big),
$$
so we can expand this determinant using (\ref{detw41D1}).  As shown before, the resulting terms in $(\lambda_1,\ldots,\lambda_r)$ are nonnegative.  Also, the resulting terms in $(x_1,\ldots,x_r)$ are of the form (\ref{vandermondederivs}) (with each $\lambda_i$ replaced by $x_i$), and hence also are nonnegative.  Moreover, it is straightforward to see that some of these terms are non-zero.  Therefore, (\ref{mixedderivdet}) is positive for $\lambda_1 > \cdots > \lambda_r$ and $x_1 > \cdots > x_r$.  

Turning to the determinant $D_2$ in (\ref{D2}), we note first that 
$$
w_4(0,x) := \lim_{\lambda \to 0+} w_4(0,x) = 1.
$$
Therefore, 
$$
D_2 = 
\left|\begin{array}{ccc}
w_4(\lambda_1,x_1) - w_4(\lambda_2,x_1) & \cdots & w_4(\lambda_1,x_r) - w_4(\lambda_2,x_r) \\
w_4(\lambda_2,x_1) - w_4(\lambda_3,x_1) & \cdots & w_4(\lambda_2,x_r) - w_4(\lambda_3,x_1)\\
\vdots & \vdots & \vdots \\
w_4(\lambda_{r-1},x_1) - w_4(\lambda_r,x_1)& \cdots & w_4(\lambda_{r-1},x_r) - w_4(\lambda_r,x_1) \\
1 & \cdots & 1
\end{array}\right|.
$$
We again apply Taylor's theorem, as in (\ref{taylor}), to each entry in rows $1,\ldots,r-1$, obtaining 
$$
D_2 = \prod_{i=1}^{r-1} (\lambda_i - \lambda_{i+1}) \cdot D_3
$$
where 
$$
D_3 = 
\left|\begin{array}{ccc}
w_{41}(\rho_1,x_1) & \cdots & w_{41}(\rho_1,x_r) \\
w_{41}(\rho_2,x_1) & \cdots & w_{41}(\rho_2,x_r) \\
\vdots & \vdots & \vdots \\
w_{41}(\rho_{r-1},x_1) & \cdots & w_{41}(\rho_{r-1},x_r) \\
1 & \cdots & 1
\end{array}\right|.
$$
Carrying out elementary column operations, subtracting column $j+1$ from column $j$, for $j=1,\ldots,r-1$, we obtain 
\begin{align*}
&D_3 \\ &= 
\left|\begin{array}{cccc}
w_{41}(\rho_1,x_1) - w_{41}(\rho_1,x_2) & \cdots & w_{41}(\rho_1,x_{r-1}) - w_{41}(\rho_1,x_r) & w_{41}(\rho_1,x_r) \\
\vdots & \vdots & \vdots & \vdots \\
w_{41}(\rho_{r-1},x_1) - w_{41}(\rho_{r-1},x_2) & \cdots & w_{41}(\rho_{r-1},x_{r-1}) - w_{41}(\rho_{r-1},x_r) & w_{41}(\rho_{r-1},x_r) \\
0 & \cdots & 0 & 1
\end{array}\right| \\
& \\
&= 
\left|\begin{array}{ccc}
w_{41}(\rho_1,x_1) - w_{41}(\rho_1,x_2) & \cdots & w_{41}(\rho_1,x_{r-1}) - w_{41}(\rho_1,x_r) \\
\vdots & \vdots & \vdots \\
w_{41}(\rho_{r-1},x_1) - w_{41}(\rho_{r-1},x_2) & \cdots & w_{41}(\rho_{r-1},x_{r-1}) - w_{41}(\rho_{r-1},x_r)
\end{array}\right|. \\
\end{align*}
Define 
%\begin{equation}
%\label{w411}
$$
w_{411}(\lambda,x) := \frac{\partial}{\partial x} w_{41}(\lambda,x) = \frac{\partial^2}{\partial \lambda \partial x} \exp\big(f(\lambda,x)\big).
$$
%\end{equation}
Applying Taylor's theorem again, we find that there exists $y_j \in (x_{j+1},x_j)$ such that 
$$
w_{41}(\rho_i,x_j) - w_{41}(\rho_i,x_{j+1}) = (x_j - x_{j+1}) w_{411}(\rho_i,y_j),
$$
$i, j = 1,\ldots,r-1$; therefore, 
$$
D_3 = \det\big((x_j - x_{j+1})  \, w_{411}(\rho_i,y_j)\big) = \prod_{j=1}^{r-1} (x_j - x_{j+1}) \cdot \det\big(w_{411}(\rho_i,y_j)\big).
$$
Noting that 
$$
\det\big(w_{411}(\rho_i,y_j)\big) = \det\left(\frac{\partial^2}{\partial \rho_i \partial y_j} \exp\big(f(\rho_i,y_j)\big)\right)
$$
is of the form (\ref{mixedderivdet}), it follows that $D_3$, and therefore $D_2$ is a product of terms, each of which is positive for $x_1 > \cdots > x_r$, $y_1 > \cdots > y_{r-1}$, and $\rho_1 > \cdots > \rho_{r-1}$; therefore $D_2 > 0$.  Consequently, $w_4(\lambda,x)$, and hence $\wtilde{w}_4(\lambda,x)$, are STP$_r$, and since $r$ was chosen arbitrarily then they both are STP$_\infty$.  
}\end{example}

\begin{example}
\label{example:w5}
{\rm 
Let 
$$
w_5(\lambda,x) = \frac{(1+\lambda)^x - 1}{\lambda x},
$$
$\lambda, x > 0$.  We show that $w_5(\lambda,x)$ is STP$_\infty$.  

We observe that 
\begin{align*}
w_5(\lambda,x) &= \lambda^{-1} \int_0^\lambda (1+t)^{x-1} \dd t \\
&= \lambda^{-1} \int_0^\infty {\bf 1}(\lambda > t( \ (1+t)^{x-1} \dd t.
\end{align*}
Recall that the function $w_2(\lambda,t) = {\bf 1}(\lambda > t)$ is TP$_\infty$; moreover, the corresponding $r \times r$ determinants are positive on an open set in $\R_+^r$.  Also, the function $w(t,x) = (1+t)^{x-1}$ is STP$_\infty$.  Therefore, by the Basic Composition Formula (\ref{bcf_continuous}), $w_5(\lambda,x)$ is STP$_\infty$.  
}\end{example}

\begin{example}
\label{example:w6}
{\rm 
Let 
$$
w_6(\lambda,x) = \frac{\lambda x}{\log(1+\lambda x)},
$$
$\lambda, x > 0$.  We shall prove that $w_6$ is STP$_\infty$ on the region $\{(\lambda,x): \lambda > 0, x > 0, \lambda x < 1\}$.  

First, we note that for all $\lambda, x > 0$, 
\begin{equation}
\label{w6integral}
w_6(\lambda,x) = \int_0^1 (1+\lambda x)^t \dd t.
\end{equation}
For $\lambda x < 1$, we expand the integrand, obtaining 
$$
(1+\lambda x)^t = \sum_{k=0}^\infty \frac{(-1)^k}{k!} (-t)_k \, \lambda^k \, x^k,
$$
where $(a)_k = a(a+1)(a+2)\cdots(a+k-1)$ is the classical rising factorial.  Substituting this series into the integral at (\ref{w6integral}) and integrating term-by-term, we obtain  
\begin{equation}
\label{w6expansion}
w_6(\lambda,x) = \sum_{k=0}^\infty \theta_k \, \lambda^k \, x^k,
\end{equation}
where 
\begin{eqnarray*}
\theta_k & = & \frac{(-1)^k }{k!} \int_0^1 (-t)_k \dd t \\
& = & \frac{1}{k!} \int_0^1 t(1-t)(2-t)\cdots(k-1-t) \dd t.
\end{eqnarray*}
This representation shows immediately that $\theta_k > 0$ for all $k \ge 0$.

By applying to (\ref{w6expansion}) the discrete version (\ref{bcf_discrete}) of the Basic Composition Formula, with $w_1(\lambda,k) = \lambda^k$, $w_2(k,x) = x^k$, and $\nu(k) = \theta_k$, we deduce that $w_6$ is STP$_r$ for all $r$; hence, $w_6(\lambda,x)$ is STP$_\infty$ on the region $\{(\lambda,x):  \lambda > 0, x > 0, \lambda x < 1\}$.  By the discrete Binet-Cauchy formula (\ref{binetcauchy_discrete}), we also obtain the representation, 
$$
\det\big(w_6(\lambda,x)\big) = \sum_{k_1 > \cdots > k_r \ge 0} \theta_{k_1} \cdots \theta_{k_r} \det(\lambda_i^{k_j}) \det(x_i^{k_j}),
$$
%Each determinant in this sum is positive for $\lambda_1 > \cdots > \lambda_r > 0$ and $x_1 > \cdots > x_r > 0$, so it follows that $w_6$ is STP$_r$ for all $r$; hence, $w_6(\lambda,x)$ is STP$_\infty$ on the region $\{(\lambda,x):  \lambda > 0, x > 0, \lambda x < 1\}$.  

As regards the total positivity properties of $w_6(\lambda,x)$ on $\R_+^2$, we calculated that for $(x_1,x_2,x_3) = (20000,0.3,0.1)$ and $(\lambda_1,\lambda_2,\lambda_3) = (3,0.4,0.1)$, 
%the $3 \times 3$ determinant, 
$$
\det\big(w_6(\lambda_i,x_j)\big)_{3 \times 3} = -5.17488 \ldots < 0.
$$
%-5.174881135976525162773482424380326027 
%Also, if $(x_1,x_2,x_3) = (3.6,1.8,0.3)$, and $(\lambda_1,\lambda_2,\lambda_3) = (1.7,1.4,0.6)$, then 
%$$
%\det\big(w_6(\lambda_i,x_j)\big)_{3 \times 3} = -0.000398258 \ldots %%48412075077605035178877304894564 < 0.
%$$
Therefore, $w_6(\lambda,x)$ is not TP$_3$ on $\R^2_+$.  
}\end{example}

\begin{example}
\label{example:w7}
{\rm
Let 
$$
w_7(\lambda,x) = \frac{\log(1+\lambda+x)}{\lambda+x}\frac{x}{\log(1+x)},
$$
$\lambda, x > 0$.  
It is straightforward to verify that 
\begin{equation}
\label{carlson}
\frac{\log(1+\lambda+x)}{\lambda+x} = \int_0^\infty (1+t)^{-1} (1+\lambda+x+t)^{-1} \dd t.
\end{equation}
We remark that this integral representation arose in work of \citet{carlson_gustafson} on the total positivity properties of mean value kernels; in the notation of Carlson and Gustafson, the function in (\ref{carlson}) is denoted by $R_{-1}(1,1;1+\lambda+x,1)$.  Writing 
$$
(1+\lambda+x+t)^{-1} = \int_0^\infty \e^{-(1+\lambda+x+t)u} \dd u,
$$
substituting this formula into (\ref{carlson}), and applying Fubini's theorem to justify an interchange of the order of integration, we obtain 
\begin{eqnarray*}
\frac{\log(1+\lambda+x)}{\lambda+x} 
& = & \int_0^\infty (1+t)^{-1} \int_0^\infty \e^{-(1+\lambda+x+t)u} \dd u \dd t \\
& = & \int_0^\infty \e^{-\lambda u} \e^{-ux} \dd\nu(u),
\end{eqnarray*}
where the positive measure $\nu$ is given explicitly by 
$$
\dd\nu(u) = \Big[\int_0^\infty \e^{-tu} (1+t)^{-1} \dd t\Big] \e^{-u} \dd u,
$$
$u > 0$.  Consequently, we have obtained an integral representation, 
$$
w_7(\lambda,x) = \frac{x}{\log(1+x)} \int_0^\infty \e^{-\lambda u} \e^{-ux} \dd\nu(u).
$$
Applying the Binet-Cauchy formula (\ref{binetcauchy_continuous}), we obtain 
\begin{multline*}
\det\big(w_7(\lambda_i,x_j)\big) \\
= \Big[\prod_{j=1}^r \frac{x_j}{\log(1+x_j)}\Big]  \idotsint\limits_{u_1 > \cdots > u_r} \det\big(\e^{-\lambda_i u_j}\big) \det\big(\e^{-u_ix_j}\big) \dd\mu(u_1) \cdots \dd\nu(u_r)
\end{multline*}
for $\lambda_1 > \cdots > \lambda_r$ and $x_1 > \cdots > x_r$.  Since the sign of each determinant in the integrand equals $(-1)^{r(r-1)/2}$ then their product is positive everywhere on the range of integration.  Therefore, $w_7(\lambda,x)$ is STP$_r$ for all $r \ge 1$, hence it is STP$_\infty$.  
}\end{example}

\section{Conclusions}
\label{sec:conclusions}

In this paper, we have explored the implications for the loading monotonicity problem of the use of higher-order totally positive weight functions for constructing weighted premiums. In doing this, we applied results from the areas of total positivity \citep{karlin} and symmetric functions \citep{macdonald}.  As a consequence, we obtained monotonicity properties of generalized weighted premiums and an upper bound under a Lipschitz hypothesis for the increase in the weighted premium in response to an increase in the loading parameter.  Further, we obtain the total positivity properties of two actuarial ratios that arise in research on combined insurance risks.  

We also examined the higher order total positivity properties of a class of kernels that have appeared in the actuarial literature.  We established the highest order of total positivity of each of these kernels, thereby adding to the collection of examples of strictly totally positive kernels.  

We related the use of weight functions that are totally positive of higher order to the degree of randomness of insured risks, and we advise insurers to relate the order of total positivity of the chosen weight function to the index of dispersion of the loss variable.  These results indicate that a broad list of TP$_\infty$, and even STP$_\infty$, weight functions is needed to develop weighted premiums for the purpose of underwriting insurable risks of any degree of randomness.

%\bigskip
%\bigskip

\section*{Acknowledgments}

We are grateful to the reviewers for helpful and constructive comments on the initial version of the manuscript.  

\iffalse
\bigskip
\bigskip

\noindent
{\bf\Large References}
\medskip
\parskip=5pt
\parindent=0pt
\fi


\begin{thebibliography}{1}
\bibliographystyle{plainnat}

%\bibpunct[; ]{(}{)}{,}{a}{}{;}

\bibitem[Brenti(1995)]{brenti}
Brenti, F. (1995).  Combinatorics and total positivity.  {\it Journal of Combinatorial Theory}, Series A, {\bf 71}, 175--218.

\bibitem[Carlson and Gustafson(1983)]{carlson_gustafson}
Carlson, B. C., and Gustafson, J. L. (1983).  Total positivity of mean values and hypergeometric functions.  {\it SIAM Journal on Mathematical Analysis}, {\bf 14}, 389--395.

\bibitem[Comtet(1974)]{comtet}
Comtet, L. (1974).  {\sl Advanced Combinatorics: The Art of Finite and Infinite Expansions}.  Reidel, Boston, MA.

\bibitem[Cox and Lewis(1966)]{cox_lewis}
Cox, D. R., and Lewis, P. A. W. (1966).  {\sl The Statistical Analysis of Series of Events}.  Methuen, London.

\bibitem[Furman and Zitikis(2008a)]{furman_zitikis08a}
Furman, E., Zitikis, R. (2008a).  Weighted premium calculation principles. {\em Insurance: Mathematics and Economics}, {\bf 42}, 459--465.

\bibitem[Furman and Zitikis(2008b)]{furman_zitikis08b}
Furman, E., and Zitikis, R. (2008b).  Monotonicity of ratios involving incomplete gamma functions with actuarial applications.  {\em Journal of Inequalities in Pure and Applied Mathematics}, Article 61, 6 pages.

\bibitem[Furman and Zitikis(2008c)]{furman_zitikis08c}
Furman, E., Zitikis, R. (2008c).  A monotonicity property of the composition of regularized and inverted-regularized gamma functions with applications. {\em Journal of Mathematical Analysis and Applications}, {\bf 348}, 971--976.

\bibitem[Furman and Zitikis(2009)]{furman_zitikis09}
Furman, E., and Zitikis, R. (2009).  Weighted pricing functionals with applications to insurance: An overview. {\em North American Actuarial Journal}, {\bf 13}, 483--496.

\bibitem[Gross and Richards(1989)]{gross_richards}
Gross, K. I., and Richards, D. St. P. (1989).  Total positivity, spherical series, and hypergeometric functions of matrix argument. {\it Journal of Approximation Theory}, {\bf 59}, 224--226.

\bibitem[Heilmann(1989)]{heilmann}
Heilmann, W.-R. (1989).  Decision theoretic foundations of credibility theory.  {\em Insurance: Mathematics and Economics}, {\bf 8}, 77--95.

\bibitem[Jones and Zitikis(2007)]{jones_zitikis}
Jones, B. L., and Zitikis, R. (2007).  Risk measures, distortion parameters, and their empirical estimation. {\em Insurance: Mathematics and Economics}, {\bf 41}, 279--297.

\bibitem[Karlin(1968)]{karlin}
Karlin, S. (1968).  {\sl Total Positivity}.  Stanford University Press, Palo Alto, CA.

\bibitem[Macdonald(1995)]{macdonald}
Macdonald, I. G. (1995).  {\sl Symmetric Functions and Hall Polynomials}, Second edition. Oxford University Press, New York.

\bibitem[Mongelli(2012)]{mongelli}
Mongelli, P. (2012).  Total positivity properties of Jacobi-Stirling numbers.  {\em Advances in Applied Mathematics}, {\bf 48}, 354--364.

\bibitem[Rinott and Saks(1993)]{rinott_saks}
 Rinott, Y., and Saks, M. (1993).  Correlation inequalities and a conjecture for permanents.  {\em Combinatorica}, {\bf 13}, 269--277.
 
\bibitem[Roman(1984)]{roman}
Roman, S. (1984).  {\sl The Umbral Calculus}. Academic Press, New York.

\bibitem[Sendov, Wang, and Zitikis(2011)]{sendov_etal}
Sendov, H. S., Wang, Y., and Zitikis, R. (2011).  Log-supermodularity of weight functions, ordering weighted losses, and the loading monotonicity of weighted premiums.  {\em Insurance: Mathematics and Economics}, {\bf 48}, 257--264.

\end{thebibliography}
\end{document}